\documentclass[12pt]{amsart}
  \usepackage{amssymb}
 \newcommand{\bee}[1]{\begin{equation}\label{#1}}
 \newcommand{\ene}{\end{equation}}

\begin{document}

MSC2010 17D99

\medskip
\begin{center}
{\Large \bf The structure of multilinear part of variety
$\widetilde{\textbf{V}}_{3}$ of Leibnitz algebras}

\medskip
{\large T. V. Skoraya}
\end{center}

\vspace{4mm} We find values of multiplicities and colength variety
of Leibniz algebras of almost polynomial growth, which is generated
by the algebra constructed with the help of the Heisenberg algebra
and its infinite-dimensional irreducible representations.

\begin{center}\textbf{1. Introduction}\end{center}

The characteristic of base field $\Phi$ throughout the paper will be
zero. This article discusses the variety
$\widetilde{\textbf{V}}_{3}$ of Leibnitz algebras and its numerical
characteristics. This variety is an analogue of the well-known
varieties $\textbf{V}_{3}$ of Lie algebras and is generated by the
Leibnitz algebra constructed with help of Heisenberg algebra and its
infinite-dimensional irreducible representations. Earlier in the
paper \cite{MSP} it is proved that the variety
$\widetilde{\textbf{V}}_{3}$ has almost polynomial growth. So the
variety $\widetilde{\textbf{V}}_{3}$ of Leibnitz algebras is extreme
in relation to the property "to have polynomial growth." The rest of
its numerical characteristics were undefined. The starting point of
the study was the result proved in the article \cite{STV-FYuYu}
about almost finite multilicities of variety $_{3}\textbf{N}$ of
Leibnitz algebras, for which $\widetilde{\textbf{V}}_{3}$ is a
subvariety. In this paper we study the multiplicities and the
colength of variety $\widetilde{\textbf{V}}_{3}.$

The first section of the paper is preliminary and contains the basic
definitions and notation needed in the sequel. In the second section
we present the structure of the generating algebra of variety
$\widetilde{\textbf{V}}_{3}$ of Leibnitz algebras and describe the
non-zero elements of a relatively free algebra of this variety,
which generate the irreducible modules of the symmetric group. In
the third section we obtain the exact values of the multiplicities.
This result made it possible to study the question of determining
the growth of other numerical characteristics, the so-called
colength, for which in the third section we have described only the
asymptotics of growth. The fourth section contains the output of the
exact formula colength of variety $\widetilde{\textbf{V}}_{3}$.

\newpage

\begin{center}\textbf{2. Basic definitions and notation.}\end{center}

A Leibnitz algebra is a vector space over a field $\Phi$ with
multiplication, which satisfies the Leibnitz identity:
$$
(xy)z\equiv (xz)y+x(yz).
$$
Probably first this class of algebras was introduced in the paper
\cite{BAM} as a generalization of concept of Lie algebra.

The determining identity of Leibnitz algebras can be represented as
follows: $x(yz)\equiv (xy)z-(xz)y.$ This form of identity allows us
any element from Leibnitz algebra to write in the form of linear
combination of elements, in which the brackets are arranged from
left to right. Therefore agree omit brackets if they are left-normed
arrangement, i.e.:
$$(((x_{1}x_{2})x_{3})\dots x_{n})=x_{1}x_{2}x_{3}\dots x_{n}.$$
For convenience we denote the operator of right multiplication for
example an element $z$ by a capital letter $Z$, assuming that
$xz=xZ.$ In particular, in our notation we obtain
$x{\underbrace{yy...y}_m}=xY^m$.

The collection of all algebras over a field $\Phi$ that satisfy a
fixed set of identities, called a variety $\textbf{V}$ of algebras
over a field $\Phi$. Note, that the system of identities can be
given implicitly. In this case the variety $\textbf{V}$  is usually
defined generating algebra given constructively.

In the paper \cite{MAI} it is proved that in the case of zero
characteristic of the base field, all information about the variety
found in the multilinear elements of its relatively free algebra.
Let $F(X,\textbf{V})$ be a relatively free algebra of variety
$\textbf{V}$ from countable set of free generators
$X=\{x_{1},x_{2},...\}$. We will denote by $P_{n}=P_{n}(\textbf{V})$
the space of all multilinear elements from generators
$x_{1},x_{2},...,x_{n}$ of algebra $F(X,\textbf{V})$. Note that for
convenience of presentation we will denote the relatively free
algebra generators also other symbols.

Let $q$ be an element of the symmetric group $S_{n}$. Assume that as
a result of the actions of the left permutation $q$ on the element
$x_{i_{1}}x_{i_{2}}\dots x_{i_{n}}$ of space $P_{n}$ we receive the
element $x_{q(i_{1})}x_{q(i_{2})}...x_{q(i_{n})}$. This sets the
action of the group $S_{n}$ on the space $P_{n}$, in consequence the
space $P_{n}$ becomes a $\Phi S_{n}$-module. This fact allows for
the study varieties of Leibnitz algebras over a field of zero
characteristic to use the well-developed theory representations of
the symmetric group.

As characteristic of the field $\Phi$ is assumed to be zero, then
the module $P_n$ is completely reducible. It is known that up to
isomorphism irreducible $\Phi S_n$-modules can be described in terms
of representations and Young diagrams.

A partition of number $n$ is a set of positive integers
$\lambda=(\lambda_1,\lambda_2,\dots, \lambda_k),$ where
$\lambda_1\ge \lambda_2\ge,\dots, \lambda_k > 0$ and
$n=\lambda_1+\lambda_2+\dots + \lambda_k.$ The partition $\lambda$
of number $n$ we will denote by $\lambda \vdash n.$ Each partition
$\lambda$ of the number $n$ of one-to-one corresponds to an Young
diagram consisting of $n$ cells in $k$ rows and containing in the
$i$-th row $\lambda_{i}$ cells.

Denoted by $\chi_{\lambda}$ the character of the irreducible
representations of the symmetric group, which corresponds to the
partition $\lambda$ of the number $n$. Then, as the module
$P_{n}(\textbf{V})$ is completely reducible, for its character is
true a decomposition:
$$\chi_{n}(\textbf{V})=\sum_{\lambda\vdash
n}m_{\lambda}(\textbf{V})\chi_{\lambda}.$$ The number of terms
$$l_{n}(\textbf{V})=\sum_{\lambda\vdash n}m_{\lambda}$$
in this sum is called a colength of the variety. Important numerical
characteristics of variety are also multiplicities
$m_{\lambda}(\textbf{V})$.

The dimension of the space $P_{n}(\bf{V})$ we denote by
$c_{n}(\bf{V})$. Let $d_{\lambda}$ be a dimension of corresponding
to $\lambda$ irreducible module. Then for introduced numerical
characteristics is carried out the following relation:
$$c_{n}(\textbf{V})=\sum_{\lambda\vdash
n}m_{\lambda}(\textbf{V})d_{\lambda}.$$ An important concept for
varieties is their growth. The growth of variety $\bf{V}$ is the
growth of the sequence of number $c_{n}(\textbf{V})$. The sequence
of number $c_{n}(\textbf{V})$ is called also a sequence of
codimensions of verbal ideal. The growth of variety is called
polynomial, if there are non-negative integers $C,m$ such that for
any $n$ is true the inequality $c_{n}(\textbf{V})<Cn^{m}$. We say
that the variety has almost polynomial growth if the growth of this
variety is not polynomial, and the growth of any of its proper
subvariety is polynomial. As previously noted the variety
$\widetilde{\textbf{V}}_{3}$ has almost polynomial growth.

Since we consider the case of zero characteristic of the base field,
then every identity is equivalent to the system of multilinear
identities, which can be obtained using the standard method of
linearization \cite{MAI}. Here is an example of this process for the
identity
$$x_{0}(xy)(xy)\equiv 0.$$
After linearization of the variable $x$ we obtain:
$$x_{0}(x_{1}y)(x_{2}y)+x_{0}(x_{2}y)(x_{1}y)\equiv 0.$$
Complete linearization is so:
$$x_{0}(x_{1}y_{1})(x_{2}y_{2})+x_{0}(x_{1}y_{2})(x_{2}y_{1})+x_{0}(x_{2}y_{1})(x_{1}y_{2})+x_{0}(x_{2}y_{2})(x_{1}y_{1})\equiv 0.$$

The space of the multilinear elements of degree $n$ of any variety
of Leibnitz algebras over a field of zero characteristic by Mashke's
theorem can be decomposed as a direct sum of irreducible submodules,
corresponding to all possible Young diagrams of $n$ cells; moreover
two modules are isomorphic if and only if they are correspond to the
same diagram. It is well-known (see, for example, \cite{GA_ZMV}),
that each from these submodules is generated by linearization of an
element $f$, which is constructed according to the Young diagram
corresponding to the partition $\lambda$ of number $n$.

To describe the structure of elements, linearization of which
generates irreducible submodules of the specified sum, we need to
introduce some notation. Recall that the standard polynomial of
degree $n$ has the form:
$$St_{n}(x_{1},x_{2},\dots ,x_{n})=\sum_{q\in
S_{n}}(-1)^{q}x_{q(1)}x_{q(2)}\dots x_{q(n)},$$ where the summation
is carried out by elements of the symmetric group, and $(-1)^{q}$ is
equal to $+1$ or $-1$ depending on the parity of permutation $q$.
Agree variables in standard polynomial denote with special symbols
above (below, wave and etc.). For example the standard polynomial of
degree $n$ in the variables $x_{1},x_{2},\dots,$ $x_{n}$ we will
write as follows:
$St_{n}=\overline{x}_{1}\overline{x}_{2}\dots\overline{x}_{n}$. It
is clear that the standard polynomial is skew symmetric. Variables
in different skew symmetric sets will be denoted by different
symbols, for example: $$\sum_{q\in S_{n}, p\in
S_{m}}(-1)^{q}(-1)^{p}x_{q(1)}x_{q(2)}\dots
x_{q(n)}y_{p(1)}y_{p(2)}\dots y_{p(m)}=$$
$$=\overline{x}_{1}\overline{x}_{2}\dots
\overline{x}_{n}\widetilde{y}_{1}\widetilde{y}_{2}\dots
\widetilde{y}_{m}.$$ Note that when the element has the same
variables in different skew symmetric sets, then its sign depends on
the parity of the permutation implicitly, therefore the variables in
this element will be called alternating. For example, element
$\overline{x}_1 \dots \overline{x}_n\widetilde{x}_1 \dots
\widetilde{x}_m$ has two alternating sets of variables.

Using the above designation we give an example of elements
corresponding, for example, to the partition $\lambda
=(m+k+l,m+k,m)$ of number $n$, where $3m+2k+l=n$, $m,k,l\geq 1$.
First, let us note that this diagram will contain three corner
cells. Recall, that the cell of the diagram is corner, if to the
right and below it there are no cells. Build a diagram corresponding
to this partition:

\vskip 0.1in
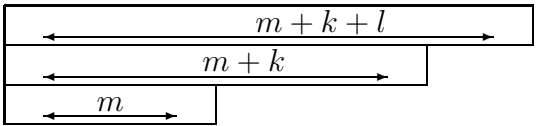
\begin{picture}(400,40)
\put(0,40){\line(1,0){200}} \put(0,25){\line(1,0){200}}
\put(0,10){\line(1,0){160}} \put(0,-5){\line(1,0){80}}
\put(200,40){\line(0,-1){15}} \put(160,25){\line(0,-1){15}}
 \put(80,10){\line(0,-1){15}}
\put(0,40){\line(0,-1){45}} \put(95,30){$m+k+l$}
\put(100,28){\vector(1,0){85}} \put(100,28){\vector(-1,0){85}}
\put(75,15){$m+k$} \put(80,13){\vector(1,0){65}}
\put(80,13){\vector(-1,0){65}} \put(35,0){$m$}
\put(40,-2){\vector(1,0){25}} \put(40,-2){\vector(-1,0){25}}
\end{picture}

\vskip 0.1in Now we construct the elements corresponding to this
diagram:
$$f_{1}=\overline{x}_{1}\overline{x}_{2}\overline{x}_{3}\widetilde{St}_{3}^{m-1}\widehat{St}_{2}^{k}X_{1}^{l},$$
$$f_{2}=\overline{x}_{1}\overline{x}_{2}\widetilde{St}_{3}^{m}\widehat{St}_{2}^{k-1}X_{1}^{l},$$
$$f_{3}=x_{1}\overline{St}_{3}^{m}\widetilde{St}_{2}^{k}X_{1}^{l-1}.$$

The described structure of elements relatively free algebra of the
variety will be used by us in the future in the proof of the
results.

\newpage \begin{center}\textbf{3. The structure of generating algebra
of variety $\widetilde{\textbf{V}}_{3}$.}\end{center}

Consider the structure of Leibnitz algebra, which generated the
variety $\widetilde{\textbf{V}}_{3}$. Let $T=\Phi[t]$ be a ring of
polynomial in the variable $t$. Consider three-dimensional
Heisenberg algebra $H$ with the besis $\{a,b,c\}$ and multiplication
$ba=-ab=c$, the product of the remaining basis elements is zero.
Well known and easy to verify that the algebra $H$ is nilpotent of
the class two Lie algebra. Transform the polynomial ring $T$ in the
right module of algebra $H$, in which the basis elements of algebra
$H$ act on the right on the polynomial $f$ from $T$ follows:
$$fa=f', fb=tf, fc=f,$$ where $f'$ is a partial derivative
of a polynomial $f$ in the variable $t$. Consider the direct sum of
vector spaces $H$ and $T$ with multiplication by the rule:
$$(x+f)(y+g)=xy+fy,$$ where $x,y$ are from $H$; $f,g$ are from $T$. Denote it by the symbol
$\widetilde{H}$. Direct verification shows that $\widetilde{H}$ is
an algebra of Leibnitz.

Thus constructed algebra $\widetilde{H}$ generates the variety
$\widetilde{\textbf{V}}_{3}$.

Determine the general form of non-zero elements of a relatively free
algebra $F(X,\widetilde{\textbf{V}}_{3})$. To do it, we will replace
the variables of these polynomials on the basis elements of algebra
$\widetilde{H}$. As a result of this replacement, we will get the
elements of the algebra $\widetilde{H}$, of which equality or
difference from zero we can check into force of the structure of
this algebra. The replacement we will choose so that  it allows to
perform a reverse replacement. This will mean that the non-zero
elements of the algebra $\widetilde{H}$ correspond to the nonzero
elements of the algebra $F(X,\widetilde{\textbf{V}}_{3})$.

Since the Heisenberg algebra $H$ is nilpotent of class two, then the
product of any its three elements is equal to zero. Consequently,
all elements of degree three and above, resulting from this
replacement, containing only the elements of the algebra $H$ are
zero. From the structure of algebra $\widetilde{H}$ follows that the
product on the left of the basis elements from $H$ by a polynomial
from $T$ is equal to zero. Therefore, the element of the algebra
$F(X,\widetilde{\textbf{V}}_{3})$ not containing in the first
alternating set the generator, which is replaced on the polynomial
from $T$, is zero. Since $T$ is regarded as a Lie algebra with zero
multiplication, then all elements of algebra $\widetilde{H}$, which
have more than one polynomial from $T$, are zero. Moreover if the
element from $F(X,\widetilde{\textbf{V}}_{3})$ has at least one
alternating set of four variables not in the first place, then as a
result of the replacement we will have element, which contains twice
one of the basic elements of algebra $H$ in the alternating set. It
is clear that such an element is zero, in account of its structure.

\newpage \begin{center}\textbf{4. The multiplicities of variety
$\widetilde{\textbf{V}}_{3}$.}\end{center}

Consider the variety $\widetilde{\textbf{V}}_{3}$ of Leibnitz
algebras and its numerical characteristics. Recall, that by
$m_{\lambda}(\widetilde{\textbf{V}}_{3})$ we denoted the
multilpicities of irreducible $\Phi S_{n}$-submodules of module
$P_{n}(\widetilde{\textbf{V}}_{3})$, which correspond to the
partition $\lambda$ of the number $n$.

\medskip\noindent
{\bf Theorem 1.} {\it Let the decomposition of character
$\chi_{n}(\widetilde{\textbf{V}}_{3})$ of the module
$P_{n}(\widetilde{\textbf{V}}_{3})$ into the integer combination of
irreducible characters $\chi_{\lambda}$ corresponding to the
partition $\lambda$ of number $n$ has the form
$$
\chi_{n}(\widetilde{\textbf{V}}_{3})=\sum_{\lambda \vdash
n}m_{\lambda}(\widetilde{\textbf{V}}_{3})\chi_{\lambda}.
$$
Then the multiplicity calculated by the formula:

$$ m_{\lambda}(\widetilde{\textbf{V}}_{3}) =
\left\{ \begin{array}{ll} {1,}& \mbox {если $\lambda = (n)$,
$\lambda =(p,p)$, $\lambda =(p,p,p)$, $\,\,n,p\ge 1,$}
\\
{ }& \mbox {или $\lambda = (p+q+r+1,p+q+1,p+1,1),$ }
\\ { }& \mbox {$\,\,p,q,r\ge 0;$}
\\
{2,}& \mbox {если $\lambda = (p+q,p)$, $\lambda = (p+q,p,p)$,}
\\ { }& \mbox {$\lambda = (p+q,p+q,p)$, $p,q\ge 1;$ }
\\
{3,}& \mbox {если $\lambda = (p+q+r,p+q,p)$,$\,\,p,q,r\geq 1.$}
\\
{0,}& \mbox {во всех остальных случаях.}
\end{array}
\right.
$$
}
\medskip\noindent
{\bf Proof.} According to the arguments given above, the space of
multilinear elements of degree $n$ any variety of Leibnitz algebras
can be decomposed as a direct sum of irreducible submodules,
corresponding to all possible Young diagrams of $n$ cells; moreover
two modules are isomorphic if and only if they are correspond to the
same diagram. In the paper \cite{ALE-MSP} it is proved that the
number of isomorphic terms in the specified sum for the space
$P_{n}(_{3}\textbf{N})$ is equal to the number of corner cells in
the corresponded Young diagram. Since the variety
$\widetilde{\textbf{V}}_{3}$ is subvariety of variety
$_{3}\textbf{N}$, then for its multilinear part the number of
isomorphic terms does not exceed the number of corner cells.

Consider the diagrams corresponded to non-zero elements, which
generates linearization irreducible $\Phi S_{n}$-submodules of the
space $P_{n}(\widetilde{\textbf{V}}_{3})$. In the paper
\cite{ALE-MSP} it is proved that these are diagrams, in which the
first column has not more than four cells and all other columns have
not more than three cells. Here are the elements corresponding to
such diagrams given in the article \cite{STV-FYuYu}.

First, consider the diagram, the first column of which has four
cells. These diagrams correspond to partitions

$\lambda=(m+1,1,1,1)$, where $m+4=n$,

$\lambda=(m+1,m+1,1,1)$, where $2m+4=n$,

$\lambda=(m+1,m+1,m+1,1)$, where $3m+4=n$,

$\lambda=(m+k+1,k+1,1,1)$, where $k\neq 1$ and $m+2k+4=n$,

$\lambda=(m+k+1,k+1,k+1,1)$, where $k\neq 1$ and $m+3k+4=n$,

$\lambda=(m+k+1,m+k+1,k+1,1)$, where $k\neq 1$ and $2m+3k+4=n$,

$\lambda=(m+k+p+1,k+p+1,p+1,1)$, where $k,p\neq 1$ and
$m+2k+3p+4=n$. Construct the corresponding elements. For partition
$\lambda=(m+1,1,1,1)$:
\begin{center} $h_{1}^{(1)}=\widetilde{x}_{1}\widetilde{x}_{2}\widetilde{x}_{3}\widetilde{x}_{4}X_{1}^{m}$,
$h_{1}^{(2)}=x_{1}\widetilde{St}_{4}X_{1}^{m-1}$;\end{center} for
partition $\lambda=(m+1,m+1,1,1)$:
\begin{center} $h_{2}^{(1)}=\widetilde{x}_{1}\widetilde{x}_{2}\widetilde{x}_{3}\widetilde{x}_{4}\overline{St}_{2}^{m}$,
$h_{2}^{(2)}=\widetilde{x}_{1}\widetilde{x}_{2}\overline{St}_{4}\widehat{St}_{2}^{m-1}$;\end{center}
for partition $\lambda=(m+1,m+1,m+1,1)$:
\begin{center} $h_{3}^{(1)}=\widetilde{x}_{1}\widetilde{x}_{2}\widetilde{x}_{3}\widetilde{x}_{4}\overline{St}_{3}^{m}$,
$h_{3}^{(2)}=\widetilde{x}_{1}\widetilde{x}_{2}\widetilde{x}_{3}\overline{St}_{4}\widehat{St}_{3}^{m-1}$;\end{center}
for partition $\lambda=(m+k+1,k+1,1,1)$:
\begin{center} $h_{4}^{(1)}=\widetilde{x}_{1}\widetilde{x}_{2}\widetilde{x}_{3}\widetilde{x}_{4}\overline{St}_{2}^{k}X_{1}^{m}$,
$h_{4}^{(2)}=\widetilde{x}_{1}\widetilde{x}_{2}\overline{St}_{4}\widehat{St}_{2}^{k-1}X_{1}^{m}$
and
$h_{4}^{(3)}=x_{1}\widetilde{St}_{4}\overline{St}_{2}^{k}X_{1}^{m-1}$;\end{center}
for partition $\lambda=(m+k+1,k+1,k+1,1)$:
\begin{center} $h_{5}^{(1)}=\widetilde{x}_{1}\widetilde{x}_{2}\widetilde{x}_{3}\widetilde{x}_{4}\overline{St}_{3}^{k}X_{1}^{m}$,
$h_{5}^{(2)}=\widetilde{x}_{1}\widetilde{x}_{2}\widetilde{x}_{3}\overline{St}_{4}\widehat{St}_{3}^{k-1}X_{1}^{m}$
and
$h_{5}^{(3)}=x_{1}\widetilde{St}_{4}\overline{St}_{3}^{k}X_{1}^{m-1}$;\end{center}
for partition $\lambda=(m+k+1,m+k+1,k+1,1)$:
\begin{center} $h_{6}^{(1)}=\widetilde{x}_{1}\widetilde{x}_{2}\widetilde{x}_{3}\widetilde{x}_{4}\overline{St}_{3}^{k}\widehat{St}_{2}^{m}$,
$h_{6}^{(2)}=\widetilde{x}_{1}\widetilde{x}_{2}\widetilde{x}_{3}\overline{St}_{4}\widehat{St}_{3}^{k-1}\widetilde{\widetilde{St}}_{2}^{m}$
and
$h_{6}^{(3)}=\widetilde{x}_{1}\widetilde{x}_{2}\overline{St}_{4}\widehat{St}_{3}^{k}\widetilde{\widetilde{St}}_{2}^{m-1}$;\end{center}
for partition $\lambda=(m+k+p+1,k+p+1,p+1,1)$:
\begin{center} $h_{7}^{(1)}=\widetilde{x}_{1}\widetilde{x}_{2}\widetilde{x}_{3}\widetilde{x}_{4}\overline{St}_{3}^{p}\widehat{St}_{2}^{k}X_{1}^{m}$,
$h_{7}^{(2)}=\widetilde{x}_{1}\widetilde{x}_{2}\widetilde{x}_{3}\overline{St}_{4}\widehat{St}_{3}^{p-1}\widetilde{\widetilde{St}}_{2}^{k}X_{1}^{m}$,
$h_{7}^{(3)}=\widetilde{x}_{1}\widetilde{x}_{2}\overline{St}_{4}\widehat{St}_{3}^{p}\widetilde{\widetilde{St}}_{2}^{k-1}X_{1}^{m}$
and
$h_{7}^{(4)}=x_{1}\widetilde{St}_{4}\overline{St}_{3}^{p}\widehat{St}_{2}^{k}X_{1}^{m}$.\end{center}

Any replacement the generators of each constructed elements, the
upper index of which is different from $(1)$, on the elements from
algebra $\widetilde{H}$ will nullify these elements according to the
arguments of the second paragraph. In the elements with the upper
index $(1)$ we will make the following replacement: $x_{1}=a$,
$x_{2}=b$, $x_{3}=c$ and $x_{4}=f$. We obtain non-zero elements.
Thus we see that to each diagram with four cells in the first column
correspond a unique irreducible submodule of the space
$P_{n}(\widetilde{\textbf{V}}_{3})$.

A similar conclusion come, having considered the diagrams of height
not more then three with one corner cell. Indeed, such diagrams
correspond to partitions $\lambda=(n)$, $\lambda=(m,m)$, where
$2m=n$ and $\lambda=(m,m,m)$, where $3m=n$. To their correspond
elements $h_{8}=x_{1}X_{1}^{n-1}$,
$h_{9}=\widetilde{x}_{1}\widetilde{x}_{2}\overline{St}_{2}^{m-1}$,
$h_{10}=\widetilde{x}_{1}\widetilde{x}_{2}\widetilde{x}_{3}\overline{St}_{3}^{m-1}$.
In this case, we can use next replacement: $x_{1}=a+f$, $x_{2}=b$,
$x_{3}=c$. So we get a nonzero elements. Thus, the linearization of
each considered element generates one irreducible submodule of the
space $P_{n}(\widetilde{\textbf{V}}_{3})$.

Now we consider diagrams of height not more than three with two
corner cells. Such diagrams correspond to partitions
$\lambda=(m+k,k)$, where $k\geq 1$ and $m+2k=n$,
$\lambda=(m+k,k,k)$, where $k\geq 1$ and $m+3k=n$ and at last
$\lambda=(m+k,m+k,k)$, where $k\geq 1$ and $2m+3k=n$. Construct
elements corresponding to these diagrams. To partition
$\lambda=(m+k,k)$ responsible elements
$h_{11}^{1}=\widetilde{x}_{1}\widetilde{x}_{2}\overline{St}_{2}^{k-1}X_{1}^{m}$
and $h_{11}^{2}=x_{1}\widetilde{St}_{2}^{k}X_{1}^{m-1}$. We show
that the elements $h_{11}^{(1)}$ and $h_{11}^{(2)}$ are linearly
independent. Assume the contrary. Suppose that there is a linear
relationship $$ \alpha_{1}h_{11}^{(1)}+\alpha_{2}h_{11}^{(2)}=0,
$$
where at least one of $\alpha_{j}$, $j=1,2$ is different from zero.
For these elements use the following replacement: $x_{1}=a$,
$x_{2}=b+f$. This substitution leads to the conclusion that the
element $h_{11}^{(2)}$ is zero, and the element $h_{11}^{(1)}$ is
different from zero. Hence, $\alpha_{1}=0$. Then it is clear that
the assumption is wrong and the elements $h_{11}^{(1)}$ and
$h_{11}^{(2)}$ are linearly independent.

To partition $\lambda=(m+k,k,k)$ correspond the elements
\begin{center}
$h_{12}^{1}=\widetilde{x}_{1}\widetilde{x}_{2}\widetilde{x}_{3}\overline{St}_{3}^{k-1}X_{1}^{m}$
and $h_{12}^{2}=x_{1}\widetilde{St}_{3}^{k}X_{1}^{m-1}$,\end{center}
and to partition $\lambda=(m+k,m+k,k)$ --- elements
\begin{center}
$h_{13}^{1}=\widetilde{x}_{1}\widetilde{x}_{2}\widetilde{x}_{3}\overline{St}_{3}^{k-1}\widehat{St}_{2}^{m}$
and
$h_{13}^{2}=\widetilde{x}_{1}\widetilde{x}_{2}\widetilde{St}_{3}^{k}\overline{St}_{2}^{m-1}$.\end{center}

Show that the elements $h_{i}^{(1)}$ and $h_{i}^{(2)}$,where $i=12,
13$, are linearly independent. Assume the contrary. Suppose that
there is a linear relationship $$
\alpha_{1}h_{i}^{(1)}+\alpha_{2}h_{i}^{(2)}=0,
$$
where at least one of $\alpha_{j}$, $j=1,2$ is different from zero.
On these elements, we introduce the following replacement:
$x_{1}=a$, $x_{2}=b$ and $x_{3}=c+f$. This exchange also resets the
elements $h_{i}^{(2)}$, and elements $h_{i}^{(1)}$ leaves non-zero.
($i=12,13$). Thus $h_{i}^{(1)}$ and $h_{i}^{(2)}$ are also linearly
independent. Therefore the linearization of each element,
corresponded to the diagram of height not more than tree with two
corner cells, generates two isomorphic irreducible submodules of
space $P_{n}(\widetilde{\textbf{V}}_{3})$.

And finally, we consider the diagrams of the last fourth type. These
include diagrams of height not more than three with three corner
cells. Such diagrams correspond to the partition
$\lambda=(m+k+p,k+p,p)$, where $k,p\geq 1$ and $m+2k+3p=n$. They
correspond to the following elements:
\begin{center} $h_{14}^{(1)}=\widetilde{x}_{1}\widetilde{x}_{2}\widetilde{x}_{3}\overline{St}_{3}^{p-1}\widehat{St}_{2}^{k}X_{1}^{m}$,
$h_{14}^{(2)}=\widetilde{x}_{1}\widetilde{x}_{2}\overline{St}_{3}^{p}\widehat{St}_{2}^{k-1}X_{1}^{m}$
and
$h_{14}^{(3)}=x_{1}\widetilde{St}_{3}^{p}\overline{St}_{2}^{k}X_{1}^{m-1}$.\end{center}
Similarly to the previous cases we have to prove their linear
independence. This was done in the paper \cite{ALE-MSP}, where for
the case of three rows was used three-dimensional Heisenberg
algebra. Therefore, the proof remains valid in the case of algebra
$\widetilde{H}$. So $h_{14}^{(1)}$, $h_{14}^{(2)}$ and
$h_{14}^{(3)}$ are linearly independent and their linearization
generates three isomorphic irreducible submodules of the space
$P_{n}(\widetilde{\textbf{V}}_{3})$. The theorem is proved.

\newpage \begin{center}\textbf{5. The asymptotic of the colength of
variety $\widetilde{\textbf{V}}_{3}$.}\end{center}

Recall that the colength of variety $\textbf{V}$ is the sum of
multiplicities $m_{\lambda}(\textbf{V})$ of this variety
$\textbf{V}$.

Since we know the multiplicity of variety
$\widetilde{\textbf{V}}_{3}$, we can now determine the nature of the
changes its colength.

\medskip\noindent
{\bf Theorem 2.} {\it For any $\varepsilon>0$ for colength variety
$\widetilde{\textbf{V}}_{3}$ of Leibnitz algebras we have the
following equality:}
$$l_{n}(\widetilde{\textbf{V}}_{3})=\frac{n^{2}}{3}+o(n^{1+\varepsilon}).$$

\medskip\noindent
{\bf Proof.} Consider the case where $n$ is large enough, for
example, $n>100$. In the proof we use Theorem 1. We consider only
the diagrams with non-zero multiplicities.

First, consider the diagram with one corner cell. Their number is
not more than three so they do not participate in the asymptotics.
The number of the diagram with two corner cell does not exceed $n$,
 and therefore is a part of the $o(n^{1+\varepsilon}).$

Estimate the number of diagrams height of three. Their number is
$\frac{n^{2}}{12}+o(n^{1+\varepsilon})$. It is a known fact.
However, describe it in detail: $\left( \begin{array}{ll} {n-2}
\\
{\ \ \ 2}
\end{array}
\right)  $ is a number of partition $n$ into three summands. Take
into account that the number of partitions that match the two or
three terms, limited to a linear function. Then the number of
partitions into three different terms is equal to
 $\left( \begin{array}{ll} {n-2}
\\
{\ \ \ 2}
\end{array}
\right)+o(n^{1+\varepsilon})$ for any $\varepsilon>0$ or
$\frac{n^{2}}{2}+o(n^{1+\varepsilon})$. The number of different
ordered partitions $3!$ times less, that is the number of diagrams
of height three with three corner cells is
$\frac{n^{2}}{12}+o(n^{1+\varepsilon})$. Given the multiplicity
received contribution to the colength
$3\cdot\frac{n^{2}}{12}+o(n^{1+\varepsilon})$.

Consider the case where the diagrams have the column of height four.
If they contain two corner of the cell, then due to the fixity of
the fourth row, their number does not exceed three.  For these same
reasons, the number of such diagrams with three corner cells is
limited by a linear function. And finally, of proved earlier, the
number of diagrams with four corner cells is equal to
$\frac{n^{2}}{12}+o(n^{1+\varepsilon})$. In this case, the
multiplicity are equal to unity and the contribution to the colength
asymptotically will be $\frac{n^{2}}{12}$. Summarizing the results,
we obtain the assertion of the theorem.

\newpage
\begin{center}\textbf{6. The colength of variety $\widetilde{\textbf{V}}_{3}$.}\end{center}

Let's go find the exact formula of colength variety
$\widetilde{\textbf{V}}_{3}$ of Leibnitz algebras. Note that theorem
2 implies that she colength can not be expressed by a polynomial.

\medskip\noindent
{\bf Theorem 3.} {\it For the colength of the variety
$\widetilde{\textbf{V}}_{3}$ of Leibnitz algebras holds following
equality:}
$$l_{n}(\widetilde{\textbf{V}}_{3})=\frac{n^{2}+n+\delta}{3},$$
where $\delta = \left\{ \begin{array}{ll} {1,}& \mbox {if $n=3k+1
$,}
\\
{0,}& \mbox {if $n\neq3k+1$.}
\end{array}
\right. $

\medskip\noindent
{\bf Proof.} Consider the sum of multiplicities
$m_{\lambda}(\widetilde{\textbf{V}}_{3})$, where $\lambda =(n)$,
$\lambda =(m+k,k)$ ($k\geq 1$) or $\lambda =(p,p)$, which correspond
to the diagrams of not more than two parts. The number of such
partitions denote by $a(n)$, and the sum of relevant multiplicities
by $l_{n}^{(2)}$. Then, if $n=2m$, are possible the partitions of
the form: $(n), (n-1,1), ..., (m,m)$. Thus we see that the number of
such diagram is $m+1$. If $n=2m+1$, then are possible follows
partitions: $(n), (n-1,1), ..., (m+1,m)$. The number of such
partitions is $m+1$. Suchwise,
$$a(n)=\left\{ \begin{array}{ll} {m+1,}& \mbox {if $n=2m$,}
\\
{m+1,}& \mbox {if $n=2m+1$,}
\end{array}
\right. $$ or $a(n)=[\frac{n}{2}]+1$.

Find the sum of the corresponding multiplicities. At the same time,
we note that the diagrams corresponding partitions $(n)$ and
$(m,m)$, have one corner cell, that is, for them
$m_{\lambda}(\widetilde{\textbf{V}}_{3})=1$. The diagrams
corresponding to the remaining partitions have two corner cells,
that is, for them $m_{\lambda}(\widetilde{\textbf{V}}_{3})=2$.
Direct calculation in both cases, we find that $l_{n}^{(2)}=n$.

Now consider the diagram of three rows. Their number we will denote
by $b(n)$, and the sum of relevant multiplicities
--- by $l_{n}^{(3)}$. Each of the considered diagrams contains at least one column of length three.
Remove these diagrams the first column. In the result of detachment
will remain diagrams with two or three rows. The number of first is
$a(n-3)$, and the  number of seconds --- $b(n-3)$. Thus, we get the
following recurrent relationship for $b(n)$:
$$b(n)=a(n-3)+b(n-3).$$ Omitting rather cumbersome calculations, we write the formula for $b(n)$:
$$b(n)=\left\{ \begin{array}{ll} {3m^{2},}& \mbox {if $n=6m$,}
\\ {3m^{2}+m,}& \mbox {if $n=6m+1$,}
\\ {3m^{2}+2m,}& \mbox {if $n=6m+2$,}
\\ {3m^{2}+3m+1,}& \mbox {if $n=6m+3$,}
\\ {3m^{2}+4m+1,}& \mbox {if $n=6m+4$,}
\\ {3m^{2}+5m+2,}& \mbox {if $n=6m+5$.}
\end{array}
\right. $$

If $n=6m$ or $n=6m+3$, then among $b(n)$ diagrams there are diagrams
with one corner cell, corresponded to the partitions $(m,m,m)$ or
$(m+1,m+1,m+1)$ relatively, for which
$m_{\lambda}(\widetilde{\textbf{V}}_{3})=1$. Also the number $b(n)$
of diagrams includes diagrams with two corner cells, for which
$m_{\lambda}(\widetilde{\textbf{V}}_{3})=2$. To find this number,
consider the partition to which they correspond. Such partitions are
divided into two types: $(m+k,k,k)\vdash n$ and $(m+k,m+k,k)\vdash
n,$ where $k\geq 1$. The kind of partitions it follows that there
can be $[\frac{n-2}{3}]$ of partitions $\lambda_{1}$, and
$[\frac{n}{6}]-1$ of partitions $\lambda_{2}$,if $n$ is divisible by
$6$, and $[\frac{n-5}{6}]+1,$ if $n$ is divisible by $6$ with the
remainder. The third and final type of diagrams, which are
considered in $b(n)$, has three corner cells and for them
$m_{\lambda}(\widetilde{\textbf{V}}_{3})=3$. Based on this
reasoning, we obtain the following formula to calculate the sum of
corresponded multiplicities:
$$l_{n}^{(3)}=\left\{ \begin{array}{ll} {9m^{2}-3m,}& \mbox {if $n=6m$,}
\\ {9m^{2},}& \mbox {if $n=6m+1$,}
\\ {9m^{2}+3m,}& \mbox {if $n=6m+2$,}
\\ {9m^{2}+6m+1,}& \mbox {if $n=6m+3$,}
\\ {9m^{2}+9m+3,}& \mbox {if $n=6m+4$,}
\\ {9m^{2}+12m+4,}& \mbox {if $n=6m+5$.}
\end{array}
\right. $$

And finally, consider the diagram of height four, number of which we
denote by $c(n)$. According to theorem 1 for them
$m_{\lambda}(\widetilde{\textbf{V}}_{3})=1$. Explain that their
number is equal to the sum $a(n-4)$ and $b(n-4)$. Remove these
diagrams the first column. Since $n$ assumes a sufficiently large,
while the number of new diagrams is the same as the previous (enough
to $n>4$). And new diagrams will be divided into two types: the
diagrams with exactly three rows, the number of which is
$$b(n-4)=\left\{
\begin{array}{ll} {3m^{2}-4m+1,}& \mbox {if $n=6m$,}
\\ {3m^{2}-3m+1,}& \mbox {if $n=6m+1$,}
\\ {3m^{2}-2m,}& \mbox {if $n=6m+2$,}
\\ {3m^{2}-m,}& \mbox {if $n=6m+3$,}
\\ {3m^{2},}& \mbox {if $n=6m+4$,}
\\ {3m^{2}+m,}& \mbox {if $n=6m+5$,}
\end{array}
\right. $$ and the diagrams with not more then three rows, which
number is
 $$a(n-4)=\left\{ \begin{array}{ll} {m-1,}& \mbox {if $n=2m$,}
\\
{m-1,}& \mbox {if $n=2m+1$.}
\end{array}
\right. $$ Then
$$с(n)=\left\{
\begin{array}{ll} {3m^{2}-m,}& \mbox {if $n=6m$,}
\\ {3m^{2},}& \mbox {if $n=6m+1$,}
\\ {3m^{2}+m,}& \mbox {if $n=6m+2$,}
\\ {3m^{2}+2m,}& \mbox {if $n=6m+3$,}
\\ {3m^{2}+3m+1,}& \mbox {if $n=6m+4$,}
\\ {3m^{2}+4m+1,}& \mbox {if $n=6m+5$,}
\end{array}
\right. $$ Because of such diagrams
$m_{\lambda}(\widetilde{\textbf{V}}_{3})=1$, the sum of
corresponding multiplicities is equal to the number of such
diagrams.

Since there are not other diagrams with non-zero multiplicities,
then summing up the results we obtain  the assertion of the theorem.
$$l_{n}(\widetilde{\textbf{V}}_{3})=\left\{ \begin{array}{ll} {12m^{2}+2m,}& \mbox {if $n=6m$,}
\\ {12m^{2}+6m+1,}& \mbox {if $n=6m+1$,}
\\ {12m^{2}+10m+2,}& \mbox {if $n=6m+2$,}
\\ {12m^{2}+14m+4,}& \mbox {if $n=6m+3$,}
\\ {12m^{2}+18m+7,}& \mbox {if $n=6m+4$,}
\\ {12m^{2}+22m+10,}& \mbox {if $n=6m+5$,}
\end{array}
\right. $$ or
$$l_{n}(\widetilde{\textbf{V}}_{3})=\frac{n^{2}+n+\delta}{3},$$
where $\delta = \left\{ \begin{array}{ll} {1,}& \mbox {if $n=3k+1
$,}
\\
{0,}& \mbox {if $n\neq3k+1$.}
\end{array}
\right. $ The theorem is proved.

\end{document}